\documentclass[12pt,a4paper]{article}
\usepackage{amsmath,amssymb,amsfonts,amsthm}

\newtheorem{theorem}{Theorem}[section]
\newtheorem{lemma}[theorem]{Lemma}
\newtheorem{proposition}[theorem]{Proposition}
\newtheorem{corollary}[theorem]{Corollary}

\newcommand{\PGammaL}{\mathop{\mathrm{P}\Gamma\mathrm{L}}}
\newcommand{\PGL}{\mathop{\mathrm{PGL}}}
\newcommand{\PSigmaL}{\mathop{\mathrm{P}\Sigma\mathrm{L}}}
\newcommand{\PSL}{\mathop{\mathrm{PSL}}}
\newcommand{\GF}{\mathop{\mathrm{GF}}}

\DeclareMathOperator{\Aut}{Aut}
\DeclareMathOperator{\Out}{Out}
\DeclareMathOperator{\soc}{soc}
\newcommand{\la}{\langle}
\newcommand{\ra}{\rangle}
\newcommand{\norml}{\vartriangleleft}

\title{Maximal subgroups of almost simple groups with socle
  $\PSL(2,q)$}
\author{Michael Giudici\\
        School Of Mathematics and Statistics\\
        The University of Western Australia\\
        35 Stirling Highway\\
        Crawley, WA 6009\\
        Australia \\
        giudici@maths.uwa.edu.au}
\date{}
\begin{document}
\maketitle
\begin{abstract}
We determine all maximal subgroups of the almost
simple groups with socle $T=\PSL(2,q)$, that is, of all groups $G$
such that $\PSL(2,q)\leqslant G\leqslant\PGammaL(2,q)$, with 
$q\geq 4$.
\end{abstract}
\section{Introduction}
The problem of determining the maximal subgroups of the almost simple
groups has a long history and has received much attention (see for
example \cite{aschbacher,kingsurvey,KL,LPS,except}). One of its most
important applications is in determining all primitive permutation
representations of such groups. 
All subgroups of $\PSL(2,q)$ were determined by Dickson in
1901 (see \cite{dickson}) and from this one can easily read off the
list of maximal subgroups. Recent combinatorial applications
\cite{johnson,edgeprim} required a list of all maximal subgroups of
any almost simple group 
with socle $\PSL(2,q)$, that is, of all groups $G$ such that
$\PSL(2,q)\leqslant G\leqslant\PGammaL(2,q)$, with $q\geq 4$. 
There are statements in the literature about maximal subgroups of
$\PGL(2,q)$ (for example \cite{huppert1,kingsurvey}) and for all such
$G$ with small values of $q$ \cite{atlas}, but despite the fact that
the general case is folklore, we are unaware of any general treatment
in the current literature. The maximal subgroups of all
low-dimensional classical groups will be determined in
\cite{colvasbook} using the framework of Aschacher's Theorem, and this
will include all groups with socle $\PSL(2,q)$. The purpose of this note
is to use Dickson's classification of the maximal subgroups of
$\PSL(2,q)$ to determine the maximal subgroups of any almost simple
group with socle $\PSL(2,q)$. Our proof follows along the lines of the
determination in \cite{kleidmanthesis} of all maximal subgroups of the
almost simple groups with socle a Suzuki group $Sz(q)$ given Suzuki's 
classification \cite{suzuki} of the maximal subgroups of $Sz(q)$. 

The group $\PGL(2,q)$ is the group of all fractional linear
transformations
$$t_{a,b,c,d}: z\mapsto \frac{az+b}{cz+d}$$
of the projective line $X=\{\infty\}\cup \GF(q)$, where $a,b,c,d\in
\GF(q)$ with $ad-bc\neq 0$. Note that
$t_{a,b,c,d}=t_{\lambda a,\lambda b,\lambda c, \lambda d}$ for all
$\lambda \in\GF(q)$. Then 
$\PSL(2,q)=\{t_{a,b,c,d}\mid ad-bc \text{ is a square}\}$. Let $\xi$ be
a primitive element of $\GF(q)$ and $\delta=t_{\xi,0,0,1}$. Then
$\PGL(2,q)=\la\PSL(2,q),\delta\ra.$ For $q=p^f$, with $p$ prime, we
can also define the map $\phi:z\mapsto z^p$  on $X$. Then
$\PGammaL(2,q)=\la \PGL(2,q),\phi\ra$ and we define 
$\PSigmaL(2,q)=\la \PSL(2,q),\phi\ra$.

Let $T$ be a finite nonabelian simple group and $G$ be an almost
simple group with socle $T$. Subgroups of $G$ containing $T$
correspond to subgroups of $G/T$ so we concentrate on the 
subgroups of $G$ which do not contain $T$. If $M$ is a maximal
subgroup of $G$ then 
$M\cap T$ is not necessarily maximal in $T$. If $M\cap T$ is not
maximal in $T$ then $M$ is called a \emph{novelty}. It is the possible
existence of novelties which requires that extra work needs to be done
to determine the maximal subgroups of $G$. The following theorem lists
all novelty maximal subgroups of almost simple groups with socle
$\PSL(2,q)$.

\begin{theorem}
\label{thm:novelties}
Let $T=\PSL(2,q)\leqslant G\leqslant \PGammaL(2,q)$ and let $M$ be a
maximal subgroup of $G$ which does not contain $T$. Then either 
$M\cap T$ is maximal in $T$, or $G$ and $M$ are given in Table
\ref{tab:novelties}.
\end{theorem}

\begin{corollary}
Let $\PGL(2,q)\leqslant G\leqslant \PGammaL(2,q)$ and suppose that
$M$ is a maximal subgroup of $G$. Then $M\cap \PGL(2,q)$ is maximal in
$\PGL(2,q)$. 
\end{corollary}

Theorem \ref{thm:novelties} is proved by combining Propositions
\ref{prn:nonlocal}, \ref{prn:r=p}, \ref{prn:rodd} and
\ref{prn:r=2}. From these propositions we can also list the maximal
subgroups of any almost simple group with socle $\PSL(2,q)$. This
provides the following two results.  

\begin{theorem}
\label{thm:PSigmaL}
Let $G=\PSigmaL(2,q)$ for $q=p^f$ for $p$ an odd prime and $f\geq 2$.
Then the maximal subgroups of $G$ which do not contain $\PSL(2,q)$ are:
\begin{enumerate}
\item the stabiliser of a point of the projective line,
\item $N_G(D_{q-1})$  for $q\neq 9$,
\item $N_G(D_{q+1})$ for $q\neq 9$,
\item $S_5$ for $p\equiv \pm 3\pmod {10}$ and $f=2$,
\item $N_G(\PSL(2,q_0))$ with $q=q_0^r$ for some prime $r$ (2
  conjugacy classes if $r=2$).
\end{enumerate}
\end{theorem}

\begin{theorem}
\label{thm:PGammaL}
Let $G=\PGammaL(2,q)$ for $q\geq 4$ not a prime. Then the maximal
subgroups of $G$ which do not contain $\PSL(2,q)$ are:
\begin{enumerate}
\item the stabiliser of a point of the projective line,
\item $N_G(D_{2(q-1)})$,
\item $N_G(D_{2(q+1)})$,
\item $N_G(\PGL(2,q_0))$ for $q=q_0^r$ with $r$ prime, $q_0\neq 2$ and
  $r$ odd if $q$ odd.
\end{enumerate}
\end{theorem}

By \cite[Theorem 2.1]{GPZ}, if $G\leqslant \PGammaL(2,q)$ is
3-transitive on the projective line then either $G$ contains
$\PGL(2,q)$, or $q=p^f$ with $p$ odd, $f$ even and 
$G=M(s,q)= \la \PSL(2,q),\phi^s\delta\ra$ for some divisor $s$ of
$f/2$. Note that $M(1,9)=M_{10}$. We can now list
the maximal subgroups of $M(s,q)$. 

\begin{theorem}
\label{thm:Msq}
Let $G=M(s,q)$ with $q=p^f$ where $p$ is an odd prime and $s$ divides
$f/2$. Then the maximal subgroups of $G$ which do not contain
$\PSL(2,q)$ are:
\begin{enumerate}
\item the stabiliser of a point of the projective line,
\item $N_G(D_{q-1})$,
\item $N_G(D_{q+1})$,
\item $N_G(\PSL(2,q_0))$ where $q=q_0^r$ with $r$ an odd
  prime. 
\end{enumerate}
\end{theorem}

\begin{table}
\begin{center}
\label{tab:novelties}
\begin{tabular}{|ll|}
\hline
$G$  &   $M$ \\
\hline\hline
$\PGL(2,7)$     &$N_G(D_6)=D_{12}$ \\
$\PGL(2,7)$     &$N_G(D_{8})=D_{16}$\\
$\PGL(2,9)$     & $N_G(D_{10})=D_{20}$\\
$\PGL(2,9)$     & $N_G(D_{8})=D_{16}$\\
$M_{10}$     & $N_G(D_{10})=C_5\rtimes C_4$\\
$M_{10}$     & $N_G(D_{8})=C_8\rtimes C_2$\\
$\PGammaL(2,9)$ &$N_G(D_{10})=C_{10}\rtimes C_4$\\
$\PGammaL(2,9)$ & $N_G(D_8)$\\
$\PGL(2,11)$    &$N_G(D_{10})=D_{20}$\\
$\PGL(2,q)$, $q=p\equiv \pm 11,19\pmod{40}$ & 
  $N_G(A_4)=S_4$\\
\hline
\end{tabular}
\end{center}
\end{table}

\section{Preliminaries}

We begin by stating Dickson's result about the maximal subgroups of
$\PSL(2,q)$. The result is divided according to the parity of $q$. 
\begin{theorem}
\label{thm:PSLqeven}
Let $q=2^f\geq 4$. Then the maximal subgroups of $\PSL(2,q)$ are:
\begin{enumerate}
\item $C_2^f\rtimes C_{q-1}$, that is, the stabiliser of a point of the
  projective line,
\item $D_{2(q-1)}$,
\item $D_{2(q+1)}$,
\item $\PGL(2,q_0)$, where $q=q_0^r$ for some prime $r$ and 
$q_0\neq 2$.
\end{enumerate}
\end{theorem}

\begin{theorem}
\label{thm:PSLqodd}
Let $q=p^f\geq 5$ with $p$ an odd prime. Then the maximal subgroups of $\PSL(2,q)$ are:
\begin{enumerate}
\item $C_p^f\rtimes C_{(q-1)/2}$, that is, the stabiliser of a point of a
  projective line,
\item $D_{q-1}$, for $q\geq 13$,
\item $D_{q+1}$, for $q\neq 7,9$,
\item $\PGL(2,q_0)$, for $q=q_0^2$ (2 conjugacy classes),
\item $\PSL(2,q_0)$, for $q=q_0^r$ where $r$ an odd prime,
\item $A_5$, for $q\equiv \pm 1\pmod {10}$, where either $q=p$ or
  $q=p^2$ and $p\equiv \pm 3\pmod {10}$ (2 conjugacy  classes),
\item $A_4$, for $q=p\equiv \pm 3\pmod 8$ and 
$q\not\equiv \pm 1 \pmod {10}$,
\item $S_4$, for $q=p\equiv \pm 1\pmod 8$ (2 conjugacy classes).
\end{enumerate}
\end{theorem}

Let $T=\PSL(2,q)$ for $q\geq 4$. Then 
$\Aut(T)=\PGammaL(2,q)=\la T,\delta,\phi\ra$. Now
$\Out(T)=\Aut(T)/T=\la \overline{\delta}\ra\times \la \overline{\phi}\ra\cong
C_{(2,q-1)}\times C_f$.  
We will frame our results in terms of the  homomorphism
$\rho:\Out(T)\rightarrow \la \overline{\delta}\ra$ defined by
$(\overline{\delta}^i,\overline{\phi}^j)\mapsto \overline{\delta}^i$
for all $i$ and $j$. If there is a unique conjugacy class of maximal
subgroups of a given isomorphism type then $\phi$ and $\delta$ fix
this class setwise. The following lemma deals with the case when there
are two conjugacy classes of a given isomorphism type.

\begin{lemma}
\label{lem:fusion}
Let $T=\PSL(2,q)\leqslant G\leqslant\PGammaL(2,q)$ with $q$ odd and
suppose that $T$ has two conjugacy classes of maximal subgroups of $T$
of the same isomorphism type. Then these two classes are fused in $G$
if and only if $\rho(G/T)\neq 1$.
\end{lemma} 
\begin{proof}
By \cite{dickson}, given the conditions on $q$ in parts (4), (6) or
(8) of Theorem \ref{thm:PSLqodd}, there is a unique conjugacy class of $S_4$
subgroups, $A_5$ subgroups and $\PGL(2,q_0)$ subgroups in $\PGL(2,q)$
and so $\delta$ fuses each pair of conjugacy classes.

The only cases where $\phi\neq 1$ and there are two classes of
isomorphic maximal subgroups are $\PGL(2,q_0)$ when $q=q_0^2$, and
$A_5$ when $q=p^2$ for $p\equiv \pm 3\pmod {10}$.  The subgroup
$\{t_{a,b,c,d}\mid a,b,c,d\in\GF(q_0)\}\cong \PGL(2,q_0)$ of $T$ is
clearly normalised by $\phi$ and so the two classes of $\PGL(2,q_0)$
subgroups are fused in $G$ if and only if $\rho(G/T)\neq 1$.
Suppose now that $q=p^2$
for $p\equiv \pm 3\pmod {10}$. Then $T$ has two conjugacy classes
$C_1,C_2$ of $A_4$ subgroups and two conjugacy classes $D_1,D_2$ of $A_5$
subgroups. Each
pair of conjugacy classes is fused in $\PGL(2,q)$. Since $A_5$ has
only one conjugacy class of $A_4$ subgroups, it follows that  
$C_1=\{R\leqslant H\mid H\in D_1, R\cong A_4\}$.
As $C_1$ and $C_2$ are fused in $\PGL(2,q)$, each $A_4$ subgroup of
$T$ is contained in an $A_5$ and so 
$C_2=\{R\leqslant H\mid H\in D_2, R\cong A_4\}$.
Hence two $A_5$ subgroups in different $T$-conjugacy classes do not
meet in an $A_4$. Let $H=A_5$ and $R\cong A_4$ be a subgroup of
$H$. Then $R$ is contained in some $\PGL(2,p)$ subgroup $S$ of $G$.
Now $S$ is centralised by some element $g\phi$ with $g\in T$ and so
$g\phi$ centralises $R$. Since the only $A_5$ subgroups containing $R$
are conjugate to $H$ it follows that $H^{g\phi}=H^{g'}$ for some  
$g'\in T$. Thus $\phi$ does not fuse $D_1$ and $D_2$ and the result follows.
\end{proof}

Given a group $H$ and prime $r$ we define $O_r(H)$ to be the largest
normal $r$-subgroup of $H$. Note that $O_r(H)$ is characteristic in
$H$. We say that $H$ is \emph{local} if it normalises an $r$-subgroup for
some prime $r$, while we say that $H$ is \emph{nonlocal} otherwise.
Let $T\norml G$ be groups and $H$ be a subgroup of $H$. We say
that $H$ \emph{extends from $T$ to $G$} if $G=TN_G(H)$.
The following lemmas are combinations of 
\cite[Lemmas 1.3.1, 1.3.2 and 1.3.3]{pomegaplus8}.

\begin{lemma}
\label{lem:facts}
Let $T$ be a nonabelian simple group and 
$T\leqslant G\leqslant \Aut(T)$. Suppose that $M$ is a maximal
subgroup of $G$ with $M$ not containing $T$ and $M_0=T\cap M$. Then
\begin{enumerate}
\item $M/M_0\cong G/T$;
\item $M=N_G(M_0)$;
\item $O_r(M_0)\norml M$;
\item If $1<K\leqslant M_0$ and $K\norml M$ then $M_0=N_T(K)$;
\item If $M_0$ is nonlocal then $C_T(\soc(M_0))=1$.
\end{enumerate}
\end{lemma}
\begin{lemma}
\label{lem:extending}
Let $T$ be a nonabelian simple group and 
$T\leqslant G\leqslant \Aut(T)$.
Then the following hold.
\begin{enumerate}
\item Suppose that $H$ is a maximal subgroup of $T$. Then $H$ extends
  from $T$ to $G$ if and only if $N_G(H)$ is a maximal subgroup of $G$.
\item A subgroup $H$ of $T$ extends to $G$ if and only if the
  $T$-conjugacy class of $H$ is the $G$-conjugacy class of $H$.
\item Suppose $1<H\leqslant K<T$, and that $K$ extends from $T$ to $G$
  and $H$ extends from $K$ to $N_G(K)$. If $H$ is
  self-normalising in $T$ then $N_G(H)\leqslant N_G(K)$.
\end{enumerate}
\end{lemma}

\section{Determining the maximal subgroups}
We first deal with the case where $M_0=M\cap T$ is nonlocal. We have the
following proposition. 
\begin{proposition}
\label{prn:nonlocal}
Let $T=\PSL(2,q)\leqslant G\leqslant\PGammaL(2,q)$ for $q\geq 4$ and
let $M$ be a maximal subgroup of $G$ not containing $T$ such that $M$
is nonlocal. Then one of the following holds.
\begin{enumerate}
\item $q=p \equiv \pm 1\pmod {10}$, $G=T$ and $M=A_5$. (2 classes)
\item $q=p^2$, $p\equiv \pm 3 \pmod {10}$, $G=T$ and $M=A_5$. (2
  classes)
\item $q=p^2$, $p\equiv \pm 3 \pmod {10}$, $G=\PSigmaL(2,q)$ and
  $M=S_5$. (2 classes)
\item $q$ even, $M=N_G(\PGL(2,q_0))$ where $q=q_0^r$ for some prime
  $r$ and $q_0\neq 2$.
\item $q$ odd, $M=N_G(\PGL(2,q_0))$, $q=q_0^2$, and 
$G\leqslant \PSigmaL(2,q)$. (2 classes) 
\item $q$ odd, $M=N_G(\PSL(2,q_0))$, $q=q_0^r$ for some odd prime $r$.
\end{enumerate}
In particular, $M$ is not a novelty. Conversely, each case listed is in
fact a maximal subgroup.
\end{proposition}
\begin{proof}
Letting $M_0=M\cap T$ and looking at the list of maximal subgroups of
$\PSL(2,q)$, we see that either $M_0=A_5$ is maximal in $T$, or
$M_0\leqslant K$, where $K=\PSL(2,q_0)$ or $\PGL(2,q_0)$ is maximal in
$T$.  Suppose first that $M_0=A_5$ is maximal in $T$. Then either
$q=p\equiv \pm 1\pmod {10}$ or $q=p^2$ and $p\equiv \pm 3\pmod{10}$. In the first instance $\Out(T)=C_2$ and the two classes of
maximal $A_5$ subgroups are fused in $\PGL(2,q)$ by Lemma
\ref{lem:fusion}. Hence we have case (1). In the second case
$\Out(T)=C_2^2$.  If $G=T$ then we are in case (2). By Lemma
\ref{lem:fusion}, if $\rho(G/T)\neq 1$ then the two classes of maximal
$A_5$ subgroups are fused in $G$ and so by Lemma \ref{lem:extending},
$N_G(A_5)$ is not maximal in $G$. This 
leaves us to consider $G=\PSigmaL(2,q)$. Since $\phi$ fixes both of
the $T$-conjugacy classes of $A_5$ subgroups it follows from Lemma
\ref{lem:extending} that $N_G(A_5)$ is maximal in $G$. Moreover, as
$A_5$ is not in a subfield group,  $\phi$ does not centralise
$A_5$. Hence $N_G(A_5)\cong S_5$ and so case (3) holds.

Suppose next that $M_0\leqslant K$, where $K=\PSL(2,q_0)$ or
$\PGL(2,q_0)$, with $K$ maximal in $T$. Then $K=C_T(a)$ for some outer
automorphism $a$ of $T$. Since $\Out(T)$ is abelian,
$[M,a]\leqslant T$. Moreover, $[M,a]$ commutes with $M_0$ and so
$[M,a]\leqslant C_T(M_0)$. By Lemma \ref{lem:facts},
$C_T(\soc(M_0))=1$ and so $[M,a]=1$. Since $M$ is maximal in $G$ it
follows that $M=C_G(a)$. Thus $M_0=C_T(a)$ and so $M_0=K$ and $M$ is
not a novelty. If $q$ is even or $q=q_0^r$ for $r$ an odd prime, then
there is a unique conjugacy class of maximal subgroups
$M_0=\PSL(2,q_0)$. Thus Lemma \ref{lem:extending} implies that
$M=N_G(M_0)$ is maximal and we get cases (4) and (6). If $q=q_0^2$
with $q$ odd then there are two classes of maximal $\PGL(2,q_0)$
subgroups in $T$ and by Lemma \ref{lem:fusion}, these are fused in $G$
if and only if $\rho(G/T)\neq 1$.  Thus $N_G(\PGL(2,q_0))$ is maximal
in $G$ if and only if $G\leqslant \PSigmaL(2,q)$. This gives case (5).
\end{proof}

When $M_0=M\cap T$ is local, there is some prime $r$ such that 
$O_r(M_0)\neq 1$. Then $O_r(M_0)$ has nontrivial centre $Z$, which is
characteristic in $M_0$. Moreover, $Z$ has a
unique maximal elementary abelian subgroup (the group generated
by all elementary abelian subgroups of $Z$) and so this is also
characteristic in $M_0$. Hence when $M_0$ is local there is an
elementary abelian $r$-subgroup $E$ of $M_0$ such that $E\norml M$.
There are three cases to consider: $r=p$, $r=2$ and $r$ is an odd
prime dividing $q\pm 1$.

\begin{proposition}
\label{prn:r=p}
Let $T=\PSL(2,q)\leqslant G\leqslant \PGammaL(2,q)$ with $q=p^f\geq 4$ for
some prime $p$, and suppose that $M$
is a maximal subgroup of $G$ which normalises an
elementary abelian $p$-subgroup $E$ of $M_0=T\cap M$. Then $M$ is the
stabiliser of a point of the projective line. In particular, $M$ is
not a novelty. Conversely, the stabiliser in $G$ of a point of the
projective line is maximal.
\end{proposition}
\begin{proof}
By Lemma \ref{lem:facts}, $M=N_G(E)$. Moreover, $E$ is contained in a
Sylow $p$-subgroup $P$ of $T$ and so is contained in some stabiliser
$K=P\rtimes C_{(q-1)/(2,q-1)}$ in $T$ of a point of the projective
line.  Since $P$ is abelian $P\leqslant M_0$, and since the only
maximal subgroup containing $P$ is $K$ we have $M_0\leqslant K$.
Moreover, as $K/P$ is abelian it follows that $M_0\norml K$. Since
$M=N_G(M_0)$ we have $M_0=K$ and so $M=N_G(K)$. Thus $M$ is the
stabiliser of a point of the projective line. Moreover, since $K$ is
maximal in $T$ it follows that $M$ is not a novelty and as $T$ has
only one conjugacy class of subgroups isomorphic to $K$, Lemma
\ref{lem:extending} implies that the stabiliser in $G$ of a point of
the projective line is maximal.  
\end{proof}

\begin{proposition}
\label{prn:rodd}
Let $T=\PSL(2,q)\leqslant G\leqslant\PGammaL(2,q)$ with $q\geq 4$ and
let $r$ be an odd prime dividing $q\pm 1$. If $M$ is a maximal
subgroup of $G$ which normalises an elementary abelian $r$-subgroup
$E$ of $M_0=M\cap T$ for some odd prime $r$ dividing $q\pm 1$ then
$M=N_G(D_{(2,q-1)(q\pm1)})$. Conversely, $N_G(D_{(2,q-1)(q\pm1)})$ is
maximal in $G$ except  
\begin{enumerate}
\item $N_G(D_4)$ when $G=\PSL(2,5)$ or $\PGL(2,5)$,
\item $D_8$ or $D_6$ when $G=\PSL(2,7)$,
\item $N_G(D_{10})$ or $N_G(D_8)$ when $G=\PSL(2,9)$ or
  $\PSigmaL(2,9)$,
\item $D_{10}$ when $G=\PSL(2,11)$.
\end{enumerate}
\end{proposition}
\begin{proof}
Since $r$ is an odd  prime dividing $q\pm 1$ it follows that
$N_T(E)=D_{(2,q-1)(q\pm 1)}$. Hence 
$M_0=D_{(2,q-1)(q\pm 1)}$ and $M=N_G(M_0)$.  For $q\neq 7,9$, Theorems
\ref{thm:PSLqeven} and \ref{thm:PSLqodd} imply that $D_{(2,q-1)(q+1)}$
is maximal in $T$ and since there is a unique conjugacy class in $T$
of such subgroups, Lemma \ref{lem:extending} implies that
$N_G(D_{(2,q-1)(q+1)})$ is maximal in $G$. Similarly, when $q\geq 13$
then $N_G(D_{(2,q-1)(q-1)})$ is maximal in $G$.  The assertions about
the maximality of $N_G(D_{(2,q-1)(q\pm 1)})$ for small values of $q$
can then be checked in \cite{atlas}.
\end{proof}

\begin{proposition}
\label{prn:r=2}
Let $T=\PSL(2,q)\leqslant G\leqslant\PGammaL(2,q)$ with $q\geq 5$ odd, and
let $M$ be a maximal subgroup of $G$ which normalises an elementary
abelian 2-subgroup $E$ of $M_0=T\cap M$.
Then one of the following holds:
\begin{enumerate}
\item $|E|=2$ and $M_0=D_{q\pm 1}$ where $2$ divides 
$\frac{q\pm 1}{2}$.
\item $|E|=4$ and $M_0=A_4$ with $q=p\equiv \pm 3\pmod 8$. Conversely,
  $N_G(A_4)$ is maximal in $G$ for these values of $q$ except when
  $G=\PSL(2,q)$ and $q\equiv \pm 1\pmod {10}$.
\item $|E|=4$, $M_0=\PSL(2,3)\cong A_4$ and $q=3^r$ for $r$ an odd
  prime. Conversely, $N_G(M_0)$ is maximal in $G$ in this case. 
\item $|E|=4$, $M_0=\PGL(2,3)\cong S_4$, $q=9$ and $G=\PSL(2,9)$ or
  $\PSigmaL(2,9)$. (2 classes) Conversely, $N_G(M_0)$ is maximal in
  $G$.
\item $|E|=4$, $M_0=S_4$, $q=p\equiv \pm 1\pmod 8$ and $G=T$ (2
  classes). Conversely $M_0$ is maximal in $G$.
\end{enumerate}
\end{proposition}
\begin{proof}
Looking at the list of maximal subgroups of $T$ we note that
$|E|=2$ or $4$. If $|E|=2$ then $M_0=N_T(E)=D_{q\pm1}$ 
where $2$ divides $\frac{q\pm 1}{2}$.  Thus we have case (1).

If $|E|=4$ then $M_0=A_4$ when $q\equiv \pm 3\pmod 8$, while $M_0=S_4$
when $q\equiv \pm 1\pmod 8$.  Suppose first that $M_0$ is contained in
a subfield group of $T$, that is $M_0\leqslant C_T(a)$ for some field
automorphism $a$. Since $\Out(T)$ is abelian, then
$[M,a]\leqslant C_T(M_0)=1$. Thus $M=C_G(a)$ and $M_0=C_T(a)$. Hence 
$q$ is a power of 3 and $C_T(a)=\PSL(2,3)\cong A_4$ or 
$\PGL(2,3)\cong S_4$. Since $M$ is maximal, $C_G(a)$ is not contained in the
centraliser of any other automorphism and so $q=3^r$ for  some prime
$r$. If $r$ is odd then $M_0=\PSL(2,3)\cong A_4$ and there is a unique
such class. Moreover, $\PSL(2,3)$ is maximal in $\PSL(2,3^r)$ for $r$
odd and so by Lemma \ref{lem:extending}, $N_G(M_0)$ is maximal in $G$.
Thus we have case (3). If $r=2$ then $M_0=\PGL(2,3)\cong S_4$ and
there are two classes of such subgroups. These classes are maximal in
$T$ and by Lemma \ref{lem:fusion} are fused in $G$ if and only if
$\rho(G/T)\neq 1$. Hence by Lemma \ref{lem:extending}, $N_G(M_0)$ is
maximal in $G$ if and only if $\rho(G/T)=1$. Thus we have case (4). 

Suppose now that $M_0$ is not contained in a subfield group.
Then by Theorem \ref{thm:PSLqodd}, either $M_0=A_4$ and 
$q=p\equiv \pm 3\pmod 8$, or $M_0=S_4$ and $q=p\equiv \pm 1\pmod 8$.
Note that the only possibilities for $G$ are then $\PSL(2,q)$ and
$\PGL(2,q)$. If $q=p\equiv \pm 1\pmod 8$ then there are two classes of
$S_4$ subgroups and these are maximal in $T$. By Lemma
\ref{lem:fusion}, they are fused in $\PGL(2,q)$ and so do not extend
to maximal subgroups of $\PGL(2,q)$. Hence we have case (5). If
$q=p\equiv \pm 3\pmod 8$ then there is a unique conjugacy class of
$A_4$ subgroups and $N_{\PGL(2,q)}(A_4)=S_4$ by \cite{dickson}. Now
$A_4$ is maximal in $T$ if and only if
$q\not\equiv \pm 1\pmod {10}$.  If $q\not\equiv \pm 1\pmod{10}$ then
Lemma \ref{lem:extending} implies that $N_{\PGL(2,q)}(A_4)$ is maximal
in $\PGL(2,q)$. If $q\equiv \pm 1\pmod{10}$ then $M_0$ is contained in
an $A_5$. However, by Lemma \ref{lem:fusion}, $\PGL(2,q)$ interchanges
the two classes of maximal $A_5$ subgroups of $T$ while there is only
one class of $A_4$ subgroups. Hence in this case we also have
$N_{\PGL(2,q)}(A_4)$ is maximal in $\PGL(2,q)$.  Thus case (2) holds.
\end{proof}

Note that the maximality of $N_G(D_{q\pm1})$ was determined in Proposition
\ref{prn:rodd}. Collating the results of Propositions
\ref{prn:nonlocal}, \ref{prn:rodd} and \ref{prn:r=2} we
can deduce Theorems \ref{thm:novelties}, \ref{thm:PSigmaL},
\ref{thm:PGammaL} and \ref{thm:Msq} follow. We also obtain the
following well known list of maximal subgroups of $\PGL(2,q)$ for $q$ odd.

\begin{theorem}
\label{thm:PGL}
Let $G=\PGL(2,q)$ with $q=p^f>3$ for some odd prime $p$. Then the
maximal subgroups of $G$ not containing $\PSL(2,q)$ are:
\begin{enumerate}
\item $C_p^f\rtimes C_{q-1}$.
\item $D_{2(q-1)}$, for $q\neq 5$.
\item $D_{2(q+1)}$.
\item $S_4$ for $q=p\equiv \pm 3\pmod 8$.
\item $\PGL(2,q_0)$ for $q=q_0^r$ with $r$ an odd prime.
\end{enumerate}
\end{theorem}

{\bf Acknowledgements}
The author thanks Colva Roney-Dougal, Derek Holt and Alice Devillers
for reading an earlier version of this paper which led to
several improvements.

\end{document}